\documentclass[11pt]{amsart}
\usepackage[T1]{fontenc}
\usepackage[utf8]{inputenc}
\usepackage[margin=1in]{geometry}
\usepackage{amsmath,amssymb,amsthm,mathtools,booktabs,array,tabularx}
\usepackage{placeins,needspace}
\usepackage{microtype,xcolor,graphicx,tikz,pgfplots}
\pgfplotsset{compat=1.18}
\usepackage[colorlinks=true,linkcolor=blue!50!black,citecolor=blue!50!black,urlcolor=blue!50!black]{hyperref}
\newtheorem{theorem}{Theorem}[section]
\newtheorem{proposition}[theorem]{Proposition}
\newtheorem{lemma}[theorem]{Lemma}
\newtheorem{corollary}[theorem]{Corollary}
\theoremstyle{definition}

\theoremstyle{remark}
\newtheorem{remark}[theorem]{Remark}
\newcommand{\R}{\mathbb R}
\newcommand{\Z}{\mathbb Z}
\newcommand{\T}{\mathbb T}
\newcommand{\E}{\mathbb E}
\newcommand{\Cov}{\operatorname{Cov}}

\newcommand{\dd}{\,\mathrm d}
\newcommand{\av}[1]{\left\langle #1\right\rangle}
\newcommand{\D}{\mathbb D}
\DeclareMathOperator{\diag}{diag}

\hypersetup{pdftitle={Fisher Geometry of Circular Wrapped Normals: Completion and Curvature},
 pdfauthor={Hyun Chul Jang, Juna Lee, Nathan Negera, Kalea Plotnick, Antony Wang},
 pdfkeywords={Fisher information, wrapped normal, metric completion, Gaussian curvature}}
\title[Circular wrapped normals: completion and curvature]{Fisher Geometry of Circular Wrapped Normals:\\ Completion and Curvature}
\author[H. C. Jang]{Hyun Chul Jang}
\address{(Hyun Chul Jang) Department of Mathematics, Texas State University, 601 University Drive, San Marcos, TX 78666-4684, USA}
\email{hcjang@txstate.edu}
\thanks{This work originated in a research project at the Honors Summer Math Camp (HSMC), organized by Mathworks at Texas State University.}
\thanks{Corresponding author: Hyun Chul Jang.}
\author[J. Lee]{Juna Lee}
\address{(Juna Lee) Princeton High School, 151 Moore Street, Princeton, NJ 08540, USA}
\email{juna10lee@gmail.com}
\author[N. Negera]{Nathan Negera}
\address{(Nathan Negera) Texas State University, 601 University Drive, San Marcos, TX 78666, USA}
\email{nathannegera@gmail.com}
\author[K. Plotnick]{Kalea Plotnick}
\address{(Kalea Plotnick) Polytechnic School, 1030 East California Boulevard, Pasadena, CA 91106, USA}
\email{kaleaplotnick@yahoo.com}
\author[A. Wang]{Antony Wang}
\address{(Antony Wang) Lambert High School, 805 Nichols Road, Suwanee, GA 30024, USA}
\email{wangantony1112@gmail.com}

\date{October 27, 2026}
\begin{document}
\raggedbottom
\begin{abstract}
We determine the intrinsic Fisher geometry of the two-parameter circular wrapped normal family. Although its finite parameters form a cylinder, its metric completion is a disk: adjoining the uniform law extends the Fisher metric to a positive-definite real-analytic metric. The Gaussian curvature vanishes at the added point, is positive nearby, and tends to $-1/2$ as the underlying standard deviation tends to zero. We prove that the curvature has exactly one nonuniform zero radius and exactly one critical radius, the latter a nondegenerate positive maximum. Opposite points sufficiently close to uniformity are joined by a unique minimizing geodesic through the uniform law, so their intrinsic distance is not attained inside the nonuniform family. In contrast, identifying the means of real normal distributions modulo a period retains a complete hyperbolic cylinder with an infinitely distant cusp. Wrapping strictly decreases Fisher information in every nonzero tangent direction, while the information deficit and its fixed-order derivatives are exponentially small as the standard deviation tends to zero. The curvature theorem is proved by convergent series and explicit analytic remainder estimates; numerical figures serve only as illustrations.
\end{abstract}
\maketitle

\section{Introduction}
Periodic data arise when an observation records an angle or a phase but not its integer lift. Starting with a real normal distribution $N(\mu,\sigma^2)$, one can impose periodicity in two different ways. Identifying $(\mu,\sigma)$ with $(\mu+kT,\sigma)$ for $k\in\Z$ gives a quotient of the Gaussian parameter space by isometries of its Fisher metric. Observing the random variable modulo $T$ instead gives a wrapped normal distribution and its own Fisher metric. Both constructions have a circular mean parameter and a positive standard deviation, but only the second discards information in the observations. We ask how this difference changes intrinsic distances, metric completion, and Gaussian curvature.

After rescaling the period to $2\pi$, write $\alpha\in\T=\R/(2\pi\Z)$ for the mean and $s>0$ for the underlying standard deviation. The Gaussian parameter quotient carries the classical metric
\[
 g_{\mathrm G}=s^{-2}(\dd\alpha^2+2\dd s^2).
\]
It is a complete hyperbolic cylinder of curvature $-1/2$. Its circles of fixed $s$ shrink as $s\to\infty$ but remain at an unbounded distance from every fixed point. In the wrapped family the same limit is the uniform law $U$ on the observation circle. Its finite radial distance follows from the Fourier representation; identifying the full intrinsic completion and the regularity of the metric at $U$ requires further analysis.

Put $q=e^{-s^2/2}$ and let $\mathcal W=\{\operatorname{Law}(X\bmod2\pi):X\sim N(\alpha,s^2)\}$, with its intrinsic Fisher distance $d_{\mathrm W}$. The first trigonometric moment $z=qe^{i\alpha}$ identifies $\mathcal W$ with the punctured disk $\D\setminus\{0\}$, where $\D=\{z\in\mathbb C:|z|<1\}$. Write the Fisher metric as $g=B(q)\dd q^2+A(q)\dd\alpha^2$, with coefficients defined in \eqref{eq:AB}.

\Needspace{9\baselineskip}
\begin{theorem}[Completion and regularity at uniformity]\label{thm:main}
The metric completion of $(\mathcal W,d_{\mathrm W})$ is obtained by adjoining the single law $U$. Under the first-moment parametrization it is isometric to $\D$ with a positive-definite real-analytic Riemannian metric extending $g$. The end $s\downarrow0$, equivalently $|z|\uparrow1$, is infinitely distant. The extended Gaussian curvature satisfies
\[
 K(U)=0,\qquad K(q)=5q^2-25q^4+O(q^6)\quad(q\downarrow0),
 \qquad \lim_{q\uparrow1}K(q)=-\frac12.
\]
\end{theorem}

The proof expresses the Fourier density as a convergent series in $z$ and $\bar z$. The resulting statistical model is regular at $z=0$, although the polar mean parameter becomes unidentifiable there. At the other end, a conditional-score identity expresses the information lost under wrapping as a positive-definite covariance. Differentiated Gaussian-tail bounds then control both radial distances and curvature. Thus the uniform limit is an ordinary interior point of the completed statistical surface. Its zero curvature is an additional property: a smooth rotationally symmetric metric need not have zero curvature at its center, and no neighborhood of $U$ is flat here.

The endpoint expansions do not determine the behavior of curvature between the two ends. The next theorem gives the global statement.

\begin{theorem}[The curvature zero and maximum]\label{thm:unique}
There is exactly one critical point $q_{\mathrm m}\in(0,1)$ of $K(q)$, and it is a nondegenerate positive global maximum. There is exactly one zero $q_*\in(0,1)$, it is simple, and $q_{\mathrm m}<q_*$. The function $K$ is strictly increasing on $(0,q_{\mathrm m})$ and strictly decreasing on $(q_{\mathrm m},1)$. It is positive on $(0,q_*)$ and negative on $(q_*,1)$.
\end{theorem}

Uniqueness here concerns the radial variable. On the completed surface the zero set consists of $U$ and one coordinate circle; the maximum is attained on one coordinate circle. The point $U$ is a strict local minimum of curvature. The proof combines Lambert series estimates near uniformity with exact Poisson summation formulas near the opposite end. The two estimates overlap and include bounds for every infinite remainder. The finite coefficient calculations are stated explicitly; no floating-point sign test or numerical root count is used.

Completion also changes which intrinsic distances are attained. For sufficiently small $q>0$, opposite first moments $qe^{i\alpha}$ and $-qe^{i\alpha}$ are joined by a unique minimizing diameter through $U$ (Proposition~\ref{prop:diameter}). Their distance is an infimum of lengths of curves in $\mathcal W$, but no curve staying in $\mathcal W$ attains it. This distinguishes the existence of a smooth completion from geodesic existence in the original statistical family.

The wrapped normal density and its Fourier representation are classical in directional statistics \cite{MardiaJupp1999}; efficient evaluation of the density is studied in \cite{Kurz2014}. Fisher information for the variance parameter of a centered wrapped normal also appears in quantum sensing \cite{Huang2024}. The Gaussian comparison belongs to classical Fisher geometry \cite{AmariNagaoka2000,CostaSantosStrapasson2015}, and the conditional-score identity is an instance of the missing-information principle \cite{Louis1982}. Fisher--Rao geometry on ambient spaces of smooth densities is studied in \cite{BruverisMichor2019}. Ambient distances alone do not identify the intrinsic completion of a subfamily, because their connecting curves need not remain within it. Our results describe the completion and global curvature of the full two-parameter circular wrapped family.

Section~\ref{sec:foundations} fixes the models and metric coefficients. Section~\ref{sec:cap} proves the analytic extension and computes its curvature. Section~\ref{sec:loss} controls the opposite end. Section~\ref{sec:comparison} identifies the completion and the minimizing diameters, completing the proof of Theorem~\ref{thm:main}. Section~\ref{sec:curvature} explains the global curvature result, whose proof is given in Appendix~\ref{app:curvature}.

\section{The two periodic constructions and their metrics}\label{sec:foundations}
\subsection{Normalization, moments, and scaling}
Let $X\sim N(\alpha,s^2)$ on $\R$ and let $\Theta=X\bmod2\pi$. For $u=\theta-\alpha$, its density is
\begin{equation}\label{eq:spatial}
 f(u;q)=\sum_{k\in\Z}\frac1{\sqrt{2\pi}s}
       \exp\left(-\frac{(u+2\pi k)^2}{2s^2}\right),
 \qquad q=e^{-s^2/2}.
\end{equation}
The notation $f(u;q)$ always denotes a density relative to angular Lebesgue measure; set $F=2\pi f$ when using the normalized angular average $\av{h}=(2\pi)^{-1}\int_0^{2\pi}h(u)\dd u$.

\begin{proposition}[Fourier representation]\label{prop:fourier}
The density in \eqref{eq:spatial} is positive, periodic, and normalized. Its moments and Fourier series are
\begin{align}
 \E[e^{in\Theta}]&=q^{n^2}e^{in\alpha},\qquad n\in\Z,\label{eq:moments}\\
 f(u;q)&=\frac1{2\pi}\sum_{n\in\Z}q^{n^2}e^{inu}
 =\frac1{2\pi}\left(1+2\sum_{n\ge1}q^{n^2}\cos(nu)\right).\label{eq:fourier}
\end{align}
Every fixed finite mixed derivative in $u,q$ converges normally on $0\le q\le\rho<1$. At $q=0$ the series defines the uniform density.
\end{proposition}
\begin{proof}
On compact positive $s$-intervals, the spatial series and its derivatives are bounded by summable Gaussian tails in $k$, uniformly on an angular period. Reindexing proves periodicity. Tonelli's theorem and $x=u+2\pi k$ give
\[
 \int_{-\pi}^{\pi}f(u;q)\dd u
 =\sum_k\int_{-\pi+2\pi k}^{\pi+2\pi k}
       \frac{e^{-x^2/(2s^2)}}{\sqrt{2\pi}s}\dd x=1.
\]
All summands are positive. Since $e^{in\Theta}=e^{inX}$, the Gaussian characteristic function gives \eqref{eq:moments}. Equivalently, integrating $e^{in\theta}$ against the spatial sum unfolds the same integral over $\R$. Absolute integrability justifies this step even for the complex integrand.

The Fourier coefficient with the convention $e^{-in\theta}$ is $(2\pi)^{-1}q^{n^2}e^{-in\alpha}$. These coefficients are absolutely summable, so their series converges uniformly to the density by uniqueness of Fourier coefficients; see \cite{SteinShakarchi2003}.

For $n^2\ge b$, a derivative of order $a$ in $u$ and $b$ in $q$ is bounded by a constant times $n^{a+2b}\rho^{n^2-b}$. Terms with $n^2<b$ have zero $b$th derivative. The resulting numerical series converges. This proves the asserted normal convergence and licenses termwise differentiation, including at $q=0$.
\end{proof}

For a fixed physical period $T>0$ and $X_T\sim N(\mu,\sigma^2)$, write
\begin{equation}\label{eq:scaling}
 \theta=\frac{2\pi}{T}(X_T\bmod T),\quad
 \alpha=\frac{2\pi\mu}{T},\quad s=\frac{2\pi\sigma}{T},\quad
 q=e^{-2\pi^2\sigma^2/T^2}.
\end{equation}
The density relative to the physical coordinate is $p_T=(2\pi/T)f$. This change of sample variable has a parameter-independent Jacobian, so it does not change Fisher information. In particular, $z=\E[e^{i\Theta}]=qe^{i\alpha}$ identifies all nonuniform members. At $q=0$ every value of $\alpha$ denotes the same uniform law.

\subsection{Fisher information and the diagonal form}
For a positive regular density $p_\eta$, define
\[
 g_{ij}(\eta)=\int\frac{(\partial_i p_\eta)(\partial_j p_\eta)}{p_\eta}
 =\E_\eta[\partial_i\log p_\eta\,\partial_j\log p_\eta].
\]
Differentiating normalization gives zero mean scores and, after a second differentiation,
\[g_{ij}=-\E[\partial_i\partial_j\log p].\]
All distances below are the intrinsic length distances of the stated metrics: connecting curves must remain in the family under consideration. Metric completions carry the topology of their completed distances.

\Needspace{11\baselineskip}
\begin{proposition}\label{prop:metric}
For $0<q<1$ the Fisher metric is
\begin{equation}\label{eq:AB}
 g=B(q)\dd q^2+A(q)\dd\alpha^2,\qquad
 A=\int_0^{2\pi}\frac{f_u^2}{f}\dd u,\quad
 B=\int_0^{2\pi}\frac{f_q^2}{f}\dd u.
\end{equation}
Both coefficients are positive and real analytic on $(0,1)$. Moreover,
\begin{equation}\label{eq:heat-q}
 qf_q=-f_{uu}.
\end{equation}
In physical coordinates the diagonal coefficients are
\[
 I_{\mu\mu}=\left(\frac{2\pi}{T}\right)^2A(q),\qquad
 I_{\sigma\sigma}=\left(\frac{4\pi^2\sigma q}{T^2}\right)^2B(q).
\]
\end{proposition}
\begin{proof}
Since $\partial_\alpha f=-f_u$, the displayed diagonal entries follow directly. The functions $f,f_q$ are even in $u$ and $f_u$ is odd, so the mixed integral vanishes on $[-\pi,\pi]$. Translation invariance removes any dependence on $\alpha$. The first Fourier harmonic shows $f_u$ and $f_q$ are not identically zero, proving positivity. On a compact subinterval of $(0,1)$, positivity supplies a uniform lower bound for $f$. The normally convergent power series in $q$ has a complex neighborhood of each real point where the denominator stays nonzero, uniformly in $u$; integration of its analytic reciprocal products proves analyticity. Multiplying the $q$ derivative of \eqref{eq:fourier} by $q$ gives exactly the negative second $u$ derivative. Finally substitute $\dd\alpha=(2\pi/T)\dd\mu$ and $\dd q=-(4\pi^2\sigma q/T^2)\dd\sigma$ into the metric.
\end{proof}

\subsection{The Gaussian parameter quotient and the comparison map}\label{subsec:quotient}
For the real family $N(\alpha,s^2)$, the two scores are $(X-\alpha)/s^2$ and $-1/s+(X-\alpha)^2/s^3$. Their covariance, computed from the centered Gaussian moments, is
\begin{equation}\label{eq:gaussianmetric}
 g_{\mathrm G}=\frac{\dd\alpha^2+2\dd s^2}{s^2}.
\end{equation}
The translations $\alpha\mapsto\alpha+2\pi k$, $k\in\Z$, act by isometries. We call the resulting cylinder with the descended metric the \emph{Gaussian parameter quotient}. Its points are orbits of real normal distributions; they do not specify a single real-line density, since $N(\alpha,s^2)$ and $N(\alpha+2\pi,s^2)$ have different densities on $\R$.

Wrapping is constant on these orbits and defines the comparison map
\[
 ([\alpha],s)\longmapsto\operatorname{Law}(X\bmod2\pi),
 \qquad X\sim N(\alpha,s^2).
\]
The first moment $e^{-s^2/2}e^{i\alpha}$ identifies both parameters, so this is a bijection onto the nonuniform wrapped family and a diffeomorphism in the displayed coordinates. The two metrics on this common parameter cylinder come from different constructions. Throughout the comparison we pull back the wrapped Fisher metric by this map; Theorem~\ref{thm:loss} proves that it is strictly smaller than \eqref{eq:gaussianmetric} in the positive-definite order. This distinction is essential when discussing completion: the map identifies the finite parameters, but need not identify Cauchy sequences for the two metrics.

\section{The uniform law is an analytic interior point}\label{sec:cap}
\begin{theorem}[Analytic extension at uniformity]\label{thm:analytic}
The map $z=qe^{i\alpha}=x+iy$ extends the wrapped family to $\D=\{z:|z|<1\}$, with $z=0$ the uniform law. The density and Fisher metric are real analytic in $(x,y)$, and the metric is positive definite everywhere. At the center,
\begin{equation}\label{eq:center}
 g_0=2(\dd x^2+\dd y^2).
\end{equation}
\end{theorem}
\begin{proof}
For $n\ge1$, set $a_n=n(n+1)/2$ and $b_n=n(n-1)/2$. Both are nonnegative integers, and
$q^{n^2}e^{in\alpha}=z^{a_n}\bar z^{b_n}$. Thus the density at a fixed observed angle is
\begin{equation}\label{eq:diskdensity}
 p_z(\theta)=\frac1{2\pi}\left[1+\sum_{n\ge1}
 \left(z^{a_n}\bar z^{b_n}e^{-in\theta}
       +\bar z^{a_n}z^{b_n}e^{in\theta}\right)\right].
\end{equation}
Temporarily replace $\bar z$ by an independent complex variable $w$. On $|z|,|w|\le\rho<1$ the $n$th summands have magnitude at most $\rho^{n^2}$, and any fixed derivative is bounded by a polynomial in $n$ times $\rho^{n^2-m}$. The series is normally convergent and holomorphic in $z,w$. Restricting to $w=\bar z$ proves real analyticity. It equals the physical positive density for $z\ne0$ and the uniform density at zero.

At any real point $z_0\in\D$, the minimum of $p_{z_0}$ over $\theta$ is positive. Uniform continuity of the complexified series then gives a small complex neighborhood in which $p$ has no zeros, uniformly in $\theta$. Thus $p_i p_j/p$ has an analytic reciprocal expansion there, and integrating proves analyticity of the metric.

If a real tangent vector $v$ has zero Fisher norm, positivity implies $\partial_vp_z(\theta)=0$ for all $\theta$. But differentiating the first moment gives
\[
 \int_0^{2\pi}e^{i\theta}\partial_vp_z(\theta)\dd\theta
 =\partial_v z=v_x+iv_y,
\]
so $v=0$. This proves positive definiteness without relying on polar coordinates. Finally,
\[
 p_z=\frac1{2\pi}(1+2x\cos\theta+2y\sin\theta)+O(|z|^4)
\]
uniformly in $\theta$. At zero the two scores are $2\cos\theta$ and $2\sin\theta$. Orthogonality gives $g_{xx}=g_{yy}=2$ and $g_{xy}=0$.
\end{proof}

\begin{remark}\label{rem:correction}
The $O(|z|^4)$ remainder in the \emph{density} does not imply the same order for the metric: its reciprocal contributes lower-order terms. In fact, the expansions below give
\begin{equation}\label{eq:cartesian}
 g=2(\dd x^2+\dd y^2)+2(x^2+y^2)(\dd x^2+\dd y^2)
       +4(x\dd x+y\dd y)^2+O(|z|^4)|\dd z|^2.
\end{equation}
Here the last term denotes a tensor whose Cartesian coefficients are $O(|z|^4)$. The quadratic terms determine the curvature at the center; the value \eqref{eq:center} alone does not determine it.
\end{remark}

\subsection{Distance and curvature at uniformity}\label{sec:series}
\begin{lemma}[Even analytic coefficients]\label{lem:even}
The coefficients $A,B$ extend as even real-analytic functions near $q=0$, and
\begin{align}
 A(q)&=2q^2+2q^4-2q^6+2q^8+O(q^{10}),\label{eq:Aseries}\\
 B(q)&=2+6q^2+2q^4+14q^6+O(q^8).\label{eq:Bseries}
\end{align}
All remainders are analytic remainders: differentiating $j$ times reduces the indicated order by $j$, for every fixed $j$ not exceeding that order.
\end{lemma}
\begin{proof}
Since $n^2\equiv n\pmod2$, $f(u;-q)=f(u+\pi;q)$. The $q$ derivative changes sign under this substitution, while the $u$ derivative does not. Squaring and integrating over one period proves $A(-q)=A(q)$ and $B(-q)=B(q)$. For complex $|q|\le\rho<1/3$,
\[
 |F(u;q)-1|\le2\sum_{n\ge1}\rho^{n^2}\le\frac{2\rho}{1-\rho}<1.
\]
Thus $F^{-1}=\sum_{j\ge0}(-h)^j$, $h=F-1$, converges normally and gives convergent power series for $A=\av{F_u^2/F}$ and $B=\av{F_q^2/F}$. This justifies both coefficient extraction and differentiated remainders.

Up to the orders needed, $h=2q\cos u+2q^4\cos2u+O(q^9)$,
$F_u=-2q\sin u-4q^4\sin2u+O(q^9)$, and
$F_q=2\cos u+8q^3\cos2u+O(q^8)$.
The omitted third harmonic first affects orders at least $q^{10}$ in $A$ and $q^8$ in $B$ even before using frequency cancellation. Taking constant Fourier coefficients in $F_u^2(-h)^j$ and $F_q^2(-h)^j$ gives Table~\ref{tab:coefficients}. Terms with $j>6$ cannot contribute to the displayed even orders: their minimal powers are $q^{j+2}$ for $A$ and $q^j$ for $B$. The possible intervening odd orders vanish by evenness. Each entry is obtained by writing $\cos nu=(e^{inu}+e^{-inu})/2$ and $\sin nu=(e^{inu}-e^{-inu})/(2i)$ and retaining exponent zero. For example, the $q^6$ term in $F_u^2(-h)$ is
\[
 -8\av{\sin^2u\cos2u}-32\av{\sin u\sin2u\cos u}
 =-8(-1/4)-32(1/4)=-6.
\]
Summing the columns proves \eqref{eq:Aseries}--\eqref{eq:Bseries}. The analytic construction above, rather than formal differentiation of an unspecified $O$ term, proves the remainder statement.
\end{proof}

\begin{table}[ht]
\centering\small
\caption{Exact contributions from the reciprocal term $(-h)^j$. All entries not shown are zero at the retained orders.}\label{tab:coefficients}
\begin{tabular}{c rrrr rrrr}\toprule
 &\multicolumn{4}{c}{$\av{F_u^2(-h)^j}$}&\multicolumn{4}{c}{$\av{F_q^2(-h)^j}$}\\
$j$&$q^2$&$q^4$&$q^6$&$q^8$&$1$&$q^2$&$q^4$&$q^6$\\\midrule
0&2&0&0&8&2&0&0&32\\
1&0&0&$-6$&0&0&0&$-18$&0\\
2&0&2&0&0&0&6&0&0\\
3&0&0&0&$-16$&0&0&0&$-88$\\
4&0&0&4&0&0&0&20&0\\
5&0&0&0&0&0&0&0&0\\
6&0&0&0&10&0&0&0&70\\\midrule
Sum&2&2&$-2$&2&2&6&2&14\\\bottomrule
\end{tabular}
\end{table}

\begin{proposition}[Exact distance to the center]\label{prop:distance}
For every $q<1$, the intrinsic distance to the uniform point is
\begin{equation}\label{eq:R}
 R(q)=\int_0^q\sqrt{B(v)}\dd v
 =\sqrt2\left(q+\frac{q^3}{2}-\frac{q^5}{8}+\frac{71q^7}{112}+O(q^9)\right).
\end{equation}
In radial arclength, the metric is $\dd R^2+G(R)^2\dd\alpha^2$, where
\begin{equation}\label{eq:G}
 G(R)=\sqrt{A(q(R))}=R-\frac{R^5}{8}+\frac{5R^7}{28}+O(R^9).
\end{equation}
For a coordinate circle obtained by fixing $q$ and varying $\alpha$, the ratio of its length to its distance from the uniform point tends to $2\pi$ as $q\downarrow0$.
\end{proposition}
\begin{proof}
The radial curve with fixed $\alpha$ has length $R(q)$. Conversely every piecewise smooth curve has speed at least $\sqrt{B(q)}|\dot q|=|\frac{\mathrm d}{\mathrm dt}R(q)|$. Integrating gives length at least the net change of $R$, so the radial curve minimizes distance from the center. For a curve passing through the center, apply this inequality away from zero and take a limit. Analytic extension from Theorem~\ref{thm:analytic} justifies that limit. Continuity of $B$ on $[0,q]$ proves finiteness.

The square-root expansion is
$\sqrt B=\sqrt2(1+\tfrac32q^2-\tfrac58q^4+\tfrac{71}{16}q^6+O(q^8))$, which integrates to \eqref{eq:R}. Put $v=R/\sqrt2$. Analytic inversion gives
\[
 q=v-\frac12v^3+\frac78v^5-\frac{295}{112}v^7+O(v^9),\qquad
 \sqrt A=\sqrt2\left(q+\frac12q^3-\frac58q^5+\frac{13}{16}q^7+O(q^9)\right).
\]
Composition yields \eqref{eq:G}. The inverse is analytic and odd because $R(q)$ is analytic and odd with positive derivative at zero. Choosing the signed analytic square root $q\sqrt{A(q)/q^2}$ shows that $G$ also extends oddly. Finally the coordinate-circle length is $2\pi G(R)$ and $G'(0)=1$.
\end{proof}

\Needspace{12\baselineskip}
\begin{proposition}[Gaussian curvature near the uniform law]\label{prop:curvature}
For $0<q<1$,
\begin{align}
 K(q)&=-\frac1{2\sqrt{AB}}\frac{\mathrm d}{\mathrm dq}
                 \left(\frac{A'}{\sqrt{AB}}\right)\label{eq:curvature}\\
 &=-\frac{A''}{2AB}+\frac{(A')^2}{4A^2B}
                         +\frac{A'B'}{4AB^2}.
\end{align}
At the center,
\begin{equation}\label{eq:Kseries}
 K(q)=5q^2-25q^4+O(q^6),\qquad
 K(R)=\frac52R^2-\frac{15}{2}R^4+O(R^6).
\end{equation}
Thus $K(0)=0$ and $K>0$ on a punctured neighborhood of the uniform law.
\end{proposition}
\begin{proof}
For $g=\dd R^2+G^2\dd\alpha^2$, the nonzero Christoffel symbols are
$\Gamma^R_{\alpha\alpha}=-GG'$ and
$\Gamma^\alpha_{R\alpha}=\Gamma^\alpha_{\alpha R}=G'/G$.
Substitution into the curvature tensor gives
$\langle\mathcal R(\partial_R,\partial_\alpha)\partial_\alpha,\partial_R\rangle=-GG''$, hence $K=-G''/G$ with the convention that a unit sphere has curvature $+1$. Since
$G'=A'/(2\sqrt{AB})$ and $\mathrm d/\mathrm dR=B^{-1/2}\mathrm d/\mathrm dq$, this is \eqref{eq:curvature}. Differentiating the bracket proves the second formula.

To see all cancellations directly from the metric coefficients,
\[
 \sqrt{AB}=2q(1+2q^2-\tfrac12q^4+4q^6+O(q^8)),\qquad
 \frac{A'}{\sqrt{AB}}=2-5q^4+10q^6+O(q^8).
\]
Consequently \eqref{eq:curvature} gives $5q^2-25q^4+O(q^6)$. The analytic remainders from Lemma~\ref{lem:even} justify the derivative and division. Alternatively differentiating \eqref{eq:G} gives $K(R)=\tfrac52R^2-\tfrac{15}{2}R^4+O(R^6)$; substituting \eqref{eq:R} gives the same answer.
\end{proof}

The equality $K(0)=0$ is stronger than smoothness at the center. For a smooth rotational metric in radial arclength with $G(R)=R+cR^3+O(R^5)$, the same curvature formula gives $K(0)=-6c$. In \eqref{eq:G} the cubic coefficient vanishes, and the first correction is $-R^5/8$. This yields zero curvature at the uniform point but positive curvature at nearby nonuniform points.

\section{Information loss and the small-standard-deviation limit}\label{sec:loss}
\subsection{Conditional scores}
\begin{theorem}[Strict missing information]\label{thm:loss}
Let $X\sim N(\alpha,s^2)$, $Y=X\bmod2\pi$, and let information matrices use the same local coordinates $(\alpha,s)$. The complete-data score is
\[
 U(X)=\begin{pmatrix}(X-\alpha)/s^2\\-1/s+(X-\alpha)^2/s^3\end{pmatrix},
 \qquad I_X=\diag(s^{-2},2s^{-2}).
\]
Then
\begin{equation}\label{eq:missing}
 \nabla\log p_Y(Y)=\E[U\mid Y],\qquad
 D(s):=I_X-I_Y=\E[\Cov(U\mid Y)]\succ0.
\end{equation}
\end{theorem}
\begin{proof}
Writing $\varphi_{\alpha,s}$ for the real Gaussian density, the conditional weight of lift $y+2\pi k$ is
$\varphi_{\alpha,s}(y+2\pi k)/p_Y(y)$. Differentiating the spatial sum is justified by Gaussian domination on a neighborhood of each $s>0$. Division by $p_Y(y)>0$ gives the conditional-score identity. The two score means are zero, and the standard Gaussian moments $\E Z^2=1$, $\E Z^3=0$, $\E Z^4=3$ give the stated $I_X$. The law of total covariance then proves the matrix equality and nonnegativity.

For strictness, fix $(a,b)\ne(0,0)$. At any observed $y$, all lift weights are positive. If the conditional variance of $aU_\alpha+bU_s$ vanished, the quadratic polynomial
\[
 P(x)=a\frac{x-\alpha}{s^2}+b\left(-\frac1s+\frac{(x-\alpha)^2}{s^3}\right)
\]
would be constant at every $x=y+2\pi k$. A polynomial constant on an infinite arithmetic progression is constant everywhere. Its quadratic and linear coefficients then give $b=0$ and $a=0$, a contradiction. Thus the conditional variance is positive for every $y$, and its expectation is positive. Finiteness follows from the Gaussian fourth moment. This proves positive definiteness of $D$.
\end{proof}

\subsection{A differentiated Gaussian-tail estimate}
Pointwise closeness of densities is insufficient to pass to curvature. The following estimate controls the differentiated information deficit and supplies the curvature and distance asymptotics as $s\downarrow0$.

\begin{lemma}[Small-standard-deviation deficit with derivatives]\label{lem:tails}
Fix $s_0>0$ and an integer $m\ge0$. There exist constants $C,N>0$, depending only on $s_0,m$, such that for every matrix entry and every $0\le j\le m$,
\begin{equation}\label{eq:tailbound}
 |\partial_s^jD_{ab}(s)|\le C s^{-N}e^{-\pi^2/(2s^2)},\qquad
 0<s\le\min(s_0,1).
\end{equation}
\end{lemma}
\begin{proof}
Translation invariance lets us set $\alpha=0$ and use $y\in[-\pi,\pi]$. Put
\[
 x_k=y+2\pi k,\quad w_k=(\sqrt{2\pi}s)^{-1}e^{-x_k^2/(2s^2)},\quad
 p=\sum_k w_k,\quad U_k=(x_k/s^2,-1/s+x_k^2/s^3)^\top.
\]
The pairwise covariance identity gives the exact formula
\begin{equation}\label{eq:pair}
 D_{ab}(s)=\frac12\int_{-\pi}^{\pi}\sum_{k\ne l}
 \frac{w_kw_l}{p}(U_{k,a}-U_{l,a})(U_{k,b}-U_{l,b})\dd y.
\end{equation}
Indeed expanding the products in the double sum with weights $w_k/p$ gives twice the conditional covariance, and multiplying by $p$ and integrating gives $D$.

We establish bounds uniform in $y$ and $s$. First, for every fixed integer $L\ge0$,
\begin{equation}\label{eq:condmom}
 \sum_k\frac{w_k}{p}(1+|x_k|)^L\le C_L.
\end{equation}
To prove this, $p\ge w_0$ and
\[
 x_k^2-y^2=4\pi ky+4\pi^2k^2
 \ge4\pi^2|k|(|k|-1).
\]
The terms $k=0,\pm1$ have bounded $|x_k|$ and weight ratio at most one; all others are bounded by a polynomial in $|k|$ times
$\exp[-2\pi^2|k|(|k|-1)/s_0^2]$, a summable sequence. This proves \eqref{eq:condmom}.

For fixed $j$, differentiating $w_k$ gives $w_k$ times a polynomial in $x_k$ and $s^{-1}$. Equation~\eqref{eq:condmom} therefore implies
$|p^{(j)}/p|\le C s^{-N_j}$. Differentiating $p\cdot p^{-1}=1$ inductively gives
\[
 |(p^{-1})^{(j)}|\le C s^{-N_j}/p.
\]
Leibniz's rule applied to each term in \eqref{eq:pair}, including the score polynomials, now bounds its $j$th derivative by
\begin{equation}\label{eq:pairbound}
 C s^{-N}(1+|k|+|l|)^L\frac{w_kw_l}{p}
 \le C s^{-N-1}(1+|k|+|l|)^L
   e^{-\max(x_k^2,x_l^2)/(2s^2)}.
\end{equation}
Here $p\ge\max(w_k,w_l)$ was used, and $N,L$ can be chosen for all $j\le m$.

If $k\ne l$, $|x_k-x_l|\ge2\pi$, so $\max(|x_k|,|x_l|)\ge\pi$. More quantitatively, let $M=\max(|k|,|l|)\ge1$. Since $|y|\le\pi$,
$\max(|x_k|,|x_l|)\ge\pi(2M-1)$. Thus
\[
 e^{-\max(x_k^2,x_l^2)/(2s^2)}
 \le e^{-\pi^2/(2s^2)}
       e^{-\pi^2((2M-1)^2-1)/(2s_0^2)}.
\]
There are at most $(2M+1)^2$ pairs at level $M$, and the last exponential makes their polynomially weighted sum finite. Integrating \eqref{eq:pairbound} over a period proves \eqref{eq:tailbound}. These same estimates on compact positive $s$-intervals dominate the differentiated double series, so every interchange of sum, derivative, and integral used in this proof is justified. Finally $s\le1$ lets us enlarge a single exponent $N$ to cover the finite collection of estimates.
\end{proof}

\Needspace{10\baselineskip}
\begin{corollary}[Metric, curvature, and distance as the standard deviation tends to zero]\label{cor:narrow}
In $(\alpha,s)$ write $g=E(s)\dd\alpha^2+C(s)\dd s^2$. Then
\begin{equation}\label{eq:narrow}
 E(s)=s^{-2}-D_{11}(s),\qquad C(s)=2s^{-2}-D_{22}(s),
\end{equation}
with the differentiated bounds of Lemma~\ref{lem:tails}. In particular,
\[
 \lim_{s\downarrow0}K(s)=-\frac12,
 \qquad \int_s^{s_1}\sqrt{C(v)}\dd v
       =\sqrt2\log(1/s)+O(1)\quad(s\downarrow0),
\]
for any fixed $s_1>0$. The end $s\downarrow0$ is at infinite intrinsic distance.
\end{corollary}
\begin{proof}
The mixed term is zero by parity, and \eqref{eq:narrow} follows from Theorem~\ref{thm:loss}. Lemma~\ref{lem:tails} through order two gives
\[
 E'= -2s^{-3}+\varepsilon_1,\quad E''=6s^{-4}+\varepsilon_2,
 \quad C'=-4s^{-3}+\varepsilon_3,
\]
where every $\varepsilon_i$ is a polynomially weighted $e^{-\pi^2/(2s^2)}$ error. Also $E,C$ are bounded below by positive multiples of $s^{-2}$ for sufficiently small $s$. Applying \eqref{eq:curvature} with $(A,B,q)$ replaced by $(E,C,s)$ therefore changes the Gaussian curvature of the unwrapped Fisher metric by at most another such error. For $E=s^{-2}$, $C=2s^{-2}$, its three terms are $-3/2$, $1/2$, and $1/2$, whose sum is $-1/2$.

Furthermore $\sqrt C=\sqrt2/s$ plus an integrable polynomially weighted exponential error. This proves the distance asymptotic. For any curve approaching $s=0$, its length is at least the change in a primitive of $\sqrt C$, so every such curve has infinite length.
\end{proof}

\section{Metric completions of the quotient and wrapped families}\label{sec:comparison}
\subsection{The Gaussian parameter quotient}
The quotient and comparison map were defined in Section~\ref{subsec:quotient}. Setting $v=\alpha/\sqrt2$ in \eqref{eq:gaussianmetric} gives $2(\dd v^2+\dd s^2)/s^2$, twice the Poincar\'e metric, with curvature $-1/2$.

Fixing $s=s_0$ and moving toward larger standard deviation along a radial geodesic gives traveled distance $\delta=\sqrt2\log(s/s_0)$. The coordinate-circle length is
\begin{equation}\label{eq:cusp}
 L_{\mathrm G}(\delta)=\frac{2\pi}{s_0}e^{-\delta/\sqrt2}.
\end{equation}
It tends to zero only as $\delta\to\infty$. This is the hyperbolic cusp. The radial lower bound on arbitrary curve lengths proves that the end is intrinsically infinitely distant, not just infinitely far along a particular curve. The opposite end $s\downarrow0$ is likewise infinitely distant. More explicitly, the length of every curve is at least $\sqrt2$ times the total variation of $\log s$. A metric Cauchy sequence therefore has $s$ bounded above and bounded away from zero. The angular mean $\alpha$ lies in a compact circle, so the sequence is contained in a compact cylinder on which the metric is smooth and positive definite. It has a convergent subsequence, and the Cauchy property gives convergence of the whole sequence. Thus the quotient is already complete and acquires no additional points.

\subsection{Completion of the wrapped family}
\begin{theorem}[Completion of the circular wrapped family]\label{thm:global}
The metric completion of the punctured-disk wrapped family is the real-analytic disk of Theorem~\ref{thm:analytic}; it is complete as a metric space. The circle $|z|=1$ is at infinite distance and contributes no completion points.
\end{theorem}
\begin{proof}
First, removing the origin does not change distances between nonzero points in the extended disk. Given a curve that enters $q<\epsilon$, replace the portion between its first entrance and last exit by an arc of $q=\epsilon$. The replacement has length at most $2\pi\sqrt{A(\epsilon)}=O(\epsilon)$. Letting $\epsilon\downarrow0$ proves equality of the punctured and extended length distances. This supplies the isometric identification needed for the completion, beyond mere density of the punctured subset.

The distance to zero is $R(q)$ by Proposition~\ref{prop:distance}. Corollary~\ref{cor:narrow} and $q=e^{-s^2/2}$ give $R(q)\to\infty$ as $q\uparrow1$. The radial length bound also shows
$|R(q_1)-R(q_2)|\le d_g(z_1,z_2)$. Thus any metric Cauchy sequence has bounded $R$ and remains in some compact Euclidean subdisk $|z|\le q_1<1$ after adjoining zero. On that subdisk the extended smooth positive-definite metric has eigenvalues bounded above and below, so it induces the Euclidean topology; a convergent subsequence and the Cauchy property prove convergence of the whole sequence. Each point of the completed disk is approached from the punctured family. This establishes the completion and completeness.

\end{proof}

\begin{remark}[What reaches the uniform law]
For every finite $s$, $q>0$ and the nonzero first moment identifies the angular mean $\alpha$. The density becomes exactly uniform only at the added limit $s=\infty$, although this limit is at finite Fisher distance. The added point is an interior point of the completed surface. The curves obtained by fixing $q$ and varying $\alpha$ are coordinate circles. In the Gaussian hyperbolic quotient they are horocycles; no such assertion is made for the wrapped surface.
\end{remark}

\begin{table}[ht]
\centering\small
\caption{The two periodic constructions in angular units.}
\begin{tabularx}{\textwidth}{@{}p{.21\textwidth}XX@{}}\toprule
 & Gaussian parameter quotient & Wrapped observation model\\\midrule
Construction & Orbits under $\alpha\mapsto\alpha+2\pi$ & Observe $X\bmod2\pi$\\
Underlying data & Real Gaussian data before taking parameter orbits & Circular data\\
Metric & $(\dd\alpha^2+2\dd s^2)/s^2$ & $B(q)\dd q^2+A(q)\dd\alpha^2$\\
End $s\to\infty$ & Infinite-distance cusp & Uniform law at finite distance\\
Completion there & No added point & One real-analytic interior point\\
Curvature & Identically $-1/2$ & $0$ at the uniform law; positive nearby; $K\to-1/2$ as $s\downarrow0$\\\bottomrule
\end{tabularx}
\end{table}

\subsection{Minimizing geodesics through the uniform law}
\begin{proposition}[Minimizing diameters and nonattainment]\label{prop:diameter}
The completed circular family is a complete Riemannian surface. For all sufficiently small $q>0$ and every $\alpha$, the unique minimizing curve, up to parametrization, from $qe^{i\alpha}$ to $-qe^{i\alpha}$ is the radial diameter through $U$. Its length is $2R(q)$. No curve confined to the nonuniform family attains this distance.
\end{proposition}
\begin{proof}
Completeness is Theorem~\ref{thm:global}; the extended metric is positive definite and real analytic by Theorem~\ref{thm:analytic}. Rotations and reflections of the observation circle are isometries. The fixed set of a reflection is a diameter, so its arclength parametrization is a geodesic: uniqueness for the geodesic initial-value problem keeps a tangent geodesic in that fixed set.

Choose a strongly convex normal neighborhood of $U$ containing $\{|z|\le\delta\}$ for some $\delta>0$; see \cite[Chapter~6]{Lee2018}. For $R(q)<R(\delta)/2$, the radial diameter has length $2R(q)$, while any connecting curve leaving that subdisk has length at least $2(R(\delta)-R(q))>2R(q)$. Within the normal neighborhood the diameter is the unique minimizing geodesic. The exit estimate makes this the unique global minimizer as well. Equality of the punctured and extended length distances was proved in Theorem~\ref{thm:global}, but the diameter contains $U$. Hence no curve in the punctured family attains that distance.
\end{proof}

In particular, the original family is not a geodesic metric space, even though its completion is a smooth complete surface. Removing the uniform point preserves distances as infima of lengths, but can remove the curves that realize those infima.

\section{Global sign and extrema of the circular curvature}\label{sec:curvature}
The endpoint estimates in Theorem~\ref{thm:main} imply positive curvature near uniformity and negative curvature near the opposite end. Theorem~\ref{thm:unique} determines the sign on the entire disk and the variation of $K(q)$ between its ends.

Its proof is given in Appendix~\ref{app:curvature}. It combines a power-series estimate near uniformity with a Poisson summation formula near the other end. The two estimates cover the full parameter interval, including analytic bounds for all remainders.

On the completed disk, $K>0$ for $0<|z|<q_*$ and $K<0$ for $q_*<|z|<1$. The zero set is the center together with the circle $|z|=q_*$. Thus a single circle separates the positive and negative curvature regions. The maximum is attained along the circle $|z|=q_{\mathrm m}$: uniqueness refers to the radius, and nondegeneracy to the radial second derivative. The center is a further critical point of curvature on the surface, a strict local minimum by $K(z)=5|z|^2+O(|z|^4)$.

Figures~\ref{fig:curvature} and~\ref{fig:cap} illustrate the curvature and the contrasting completions. The decimal values are numerical approximations used only in the figures.

\begin{figure}[ht]
\centering
\begin{tikzpicture}
\begin{axis}[width=.92\textwidth,height=6cm,xlabel={Angular standard deviation $s$},ylabel={$K$},xmin=.3,xmax=3.2,ymin=-.57,ymax=.46,grid=major,grid style={gray!15},legend style={font=\footnotesize,at={(.98,.03)},anchor=south east}]
\addplot[blue!55!black,thick] coordinates {(0.3,-0.49999999999999933) (0.31620111731843575,-0.5000000000000013) (0.3324022346368715,-0.49999999999999756) (0.34860335195530723,-0.4999999999999223) (0.364804469273743,-0.49999999999788636) (0.38100558659217876,-0.4999999999630923) (0.3972067039106145,-0.4999999995510146) (0.4134078212290503,-0.4999999959747724) (0.42960893854748605,-0.49999997216234116) (0.4458100558659218,-0.49999984581615764) (0.4620111731843575,-0.4999992945468087) (0.4782122905027933,-0.4999972639471967) (0.49441340782122906,-0.499990808353862) (0.5106145251396648,-0.49997276146489544) (0.5268156424581005,-0.49992769015225247) (0.5430167597765363,-0.4998257622729111) (0.5592178770949721,-0.49961460783694844) (0.5754189944134078,-0.49920990533114007) (0.5916201117318436,-0.4984861006014769) (0.6078212290502794,-0.4972691353339579) (0.6240223463687151,-0.4953331405437106) (0.6402234636871509,-0.4924026683863997) (0.6564245810055866,-0.48816126159316364) (0.6726256983240224,-0.48226616970665237) (0.6888268156424582,-0.47436803640363645) (0.7050279329608938,-0.4641336042790388) (0.7212290502793296,-0.45126904750723085) (0.7374301675977654,-0.43554149689200916) (0.7536312849162011,-0.4167966324256971) (0.7698324022346369,-0.39497079274718383) (0.7860335195530725,-0.37009676622734566) (0.8022346368715083,-0.3423031609080499) (0.8184357541899441,-0.31180789673807774) (0.8346368715083798,-0.2789068532454787) (0.8508379888268156,-0.24395900567200401) (0.8670391061452514,-0.20736949286209794) (0.8832402234636871,-0.16957200641449344) (0.8994413407821229,-0.1310117136683866) (0.9156424581005587,-0.0921296735117555) (0.9318435754189944,-0.05334941799916981) (0.9480446927374302,-0.015066091029326989) (0.964245810055866,0.022361715641985525) (0.9804469273743017,0.05861749306735142) (0.9966480446927375,0.09342992753762769) (1.0128491620111733,0.12657390851120998) (1.029050279329609,0.15786979777659083) (1.0452513966480448,0.18718133948071847) (1.0614525139664805,0.2144125759755574) (1.0776536312849163,0.23950410037224823) (1.0938547486033519,0.2624289312896194) (1.1100558659217876,0.2831882450416233) (1.1262569832402234,0.30180715027675875) (1.1424581005586592,0.31833064318320414) (1.158659217877095,0.3328198398192357) (1.1748603351955307,0.3453485468997006) (1.1910614525139664,0.3560002036868585) (1.2072625698324022,0.36486520519666044) (1.223463687150838,0.37203860012920087) (1.2396648044692737,0.3776181449840118) (1.2558659217877095,0.3817026879054549) (1.2720670391061453,0.3843908511177827) (1.288268156424581,0.38577997862582725) (1.3044692737430168,0.3859653155409981) (1.3206703910614526,0.38503938641294877) (1.3368715083798883,0.3830915418711255) (1.353072625698324,0.3802076453719758) (1.3692737430167599,0.3764698746415237) (1.3854748603351956,0.3719566153152344) (1.4016759776536314,0.36674242716121286) (1.4178770949720672,0.3608980660417591) (1.434078212290503,0.35449054735773045) (1.4502793296089387,0.34758323909440303) (1.4664804469273744,0.3402359747317678) (1.4826815642458102,0.3325051781871029) (1.498882681564246,0.32444399463269935) (1.5150837988826815,0.31610242248223486) (1.5312849162011173,0.30752744308165614) (1.547486033519553,0.2987631456893223) (1.5636871508379888,0.2898508462022118) (1.5798882681564246,0.28082919879390666) (1.5960893854748603,0.27173430020041023) (1.612290502793296,0.2625997868256851) (1.6284916201117319,0.2534569251665679) (1.6446927374301676,0.2443346962840356) (1.6608938547486034,0.23525987519345049) (1.6770949720670392,0.226257106117147) (1.693296089385475,0.2173489745599635) (1.7094972067039107,0.20855607713246216) (1.7256983240223465,0.19989708998061362) (1.7418994413407822,0.19138883658142802) (1.758100558659218,0.1830463555494306) (1.7743016759776538,0.17488296897511457) (1.7905027932960895,0.1669103516809909) (1.8067039106145253,0.15913860165975824) (1.822905027932961,0.15157631182409492) (1.8391061452513968,0.14423064309585176) (1.8553072625698326,0.13710739875062417) (1.8715083798882681,0.1302110998489836) (1.887709497206704,0.12354506151170086) (1.9039106145251397,0.11711146973944775) (1.9201117318435754,0.1109114584302695) (1.9363128491620112,0.10494518622433224) (1.952513966480447,0.09921191278547181) (1.9687150837988827,0.09371007413660748) (1.9849162011173185,0.08843735665242902) (2.001117318435754,0.08339076936277756) (2.01731843575419,0.07856671421028705) (2.0335195530726256,0.07396105396418462) (2.0497206703910615,0.06956917751085356) (2.065921787709497,0.06538606228709609) (2.082122905027933,0.061406333662251456) (2.0983240223463686,0.05762432111319149) (2.1145251396648046,0.05403411107503642) (2.13072625698324,0.05062959639695408) (2.146927374301676,0.04740452235633885) (2.1631284916201117,0.04435252923442912) (2.1793296089385477,0.04146719146829603) (2.195530726256983,0.03874205344594839) (2.2117318435754187,0.03617066200852137) (2.2279329608938547,0.033746595772014) (2.2441340782122903,0.031463491375376185) (2.2603351955307263,0.02931506679355933) (2.276536312849162,0.027295141849931204) (2.292737430167598,0.02539765609920508) (2.3089385474860333,0.02361668421521168) (2.3251396648044693,0.021946449046821526) (2.341340782122905,0.020381332534266328) (2.357541899441341,0.01891588461367011) (2.3737430167597764,0.017544830271406964) (2.389944134078212,0.016263074940720212) (2.406145251396648,0.01506570834564469) (2.4223463687150835,0.013948006976125527) (2.4385474860335195,0.012905435298293133) (2.454748603351955,0.011933645869077747) (2.470949720670391,0.011028478459524704) (2.4871508379888265,0.010185958303479348) (2.5033519553072625,0.009402293594015898) (2.519553072625698,0.008673872312630504) (2.535754189944134,0.007997258513302796) (2.5519553072625696,0.007369188116693071) (2.5681564245810056,0.006786564307619702) (2.584357541899441,0.006246452641977873) (2.600558659217877,0.00574607584074549) (2.6167597765363126,0.005282808453704657) (2.6329608938547486,0.004854171332963064) (2.649162011173184,0.004457826039244495) (2.66536312849162,0.004091569188509059) (2.6815642458100557,0.003753326766446019) (2.6977653631284917,0.0034411485153817978) (2.7139664804469272,0.0031532023090699024) (2.7301675977653628,0.002887768675160274) (2.7463687150837988,0.0026432352732352893) (2.7625698324022343,0.0024180917046661143) (2.7787709497206703,0.002210924268031933) (2.794972067039106,0.002020410879659318) (2.811173184357542,0.0018453162513283328) (2.8273743016759774,0.001684487150674041) (2.8435754189944134,0.0015368477968960015) (2.859776536312849,0.0014013954759652947) (2.875977653631285,0.001277196352845289) (2.8921787709497204,0.0011633813953515482) (2.9083798882681564,0.001059142526033563) (2.924581005586592,0.000963728860805757) (2.940782122905028,0.0008764432501739972) (2.9569832402234635,0.0007966388834832472) (2.9731843575418995,0.000723716071437952) (2.989385474860335,0.0006571192125741643) (3.005586592178771,0.0005963340887495949) (3.0217877094972065,0.0005408848446850634) (3.0379888268156425,0.0004903315749251735) (3.054189944134078,0.000444267932080189) (3.0703910614525136,0.0004023187637516301) (3.0865921787709496,0.0003641380409591388) (3.102793296089385,0.0003294065444521266) (3.118994413407821,0.0002978304462677439) (3.1351955307262567,0.0002691390592579612) (3.1513966480446927,0.00024308346373547174) (3.167597765363128,0.0002194343863775039) (3.183798882681564,0.00019798185174102478) (3.2,0.00017853225928576855)};
\addlegendentry{Wrapped normal}
\addplot[red!70!black,dashed,domain=.3:3.2] {-.5};
\addlegendentry{Gaussian parameter quotient}
\addplot[gray,thin,domain=.3:3.2] {0};
\end{axis}
\end{tikzpicture}
\caption{Numerical illustration of the wrapped Gaussian curvature and the constant quotient curvature. The wrapped zero occurs at $s\approx0.9545$. Theorem~\ref{thm:unique} proves uniqueness of the zero and maximum.}\label{fig:curvature}
\end{figure}
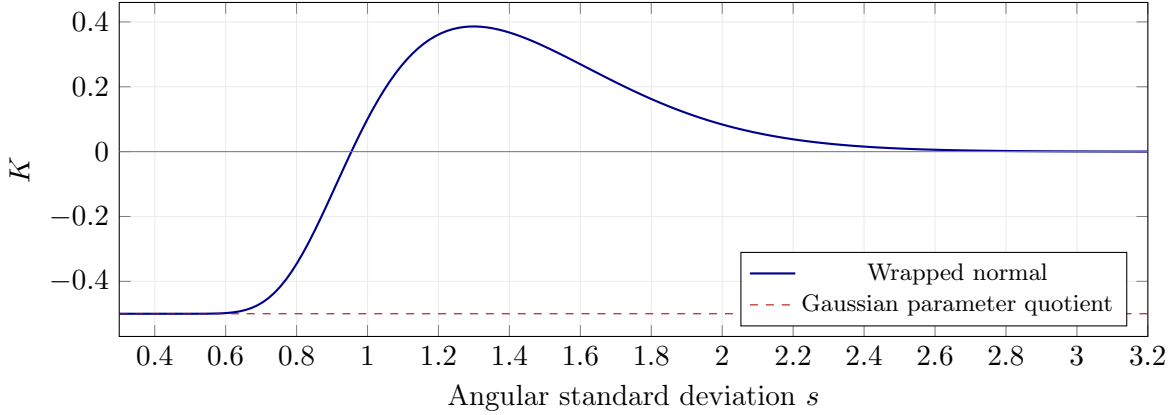
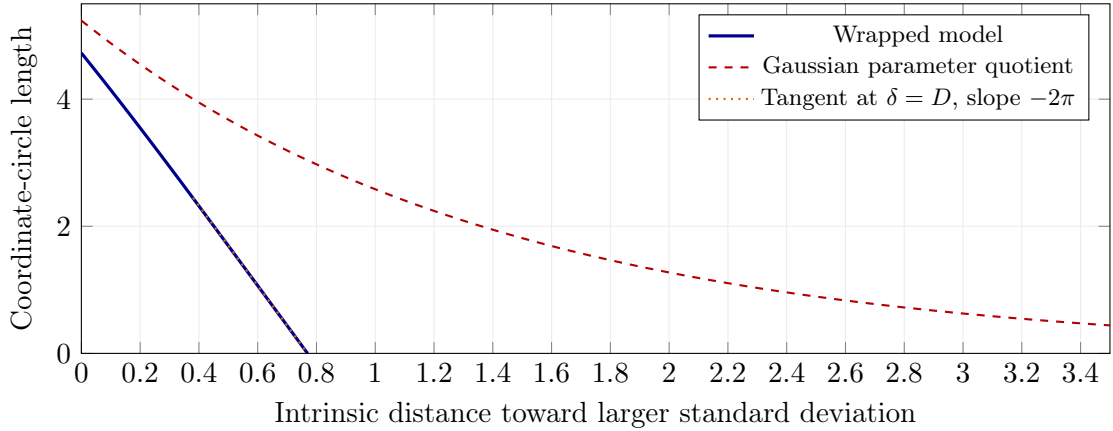
\begin{figure}[ht]
\centering
\begin{tikzpicture}
\begin{axis}[width=.92\textwidth,height=6.2cm,xlabel={Intrinsic distance toward larger standard deviation},ylabel={Coordinate-circle length},xmin=0,xmax=3.5,ymin=0,ymax=5.5,grid=major,grid style={gray!15},legend style={font=\footnotesize,at={(.98,.97)},anchor=north east}]
\addplot[blue!55!black,very thick] coordinates {(0.0,4.725738577373219) (0.015589324095328316,4.636254168090057) (0.03102426731112584,4.54724753333375) (0.04630951129164884,4.458705259853601) (0.06144952787908886,4.370615850107153) (0.07644859434557327,4.28296947230111) (0.09131080751990617,4.195757734200854) (0.1060400968929831,4.108973478715422) (0.12064023678807134,4.0226105993984085) (0.13511485767430675,3.9366638741351943) (0.14946745669468853,3.851128815411928) (0.1637014074734371,3.766001535681825) (0.1778199692617638,3.6812786264588477) (0.1918262954757911,3.59695704987765) (0.20572344167554057,3.5130340415616255) (0.2195143730294662,3.429507023737956) (0.2332019713049741,3.3463735276297815) (0.24678904142164337,3.263631124241024) (0.26027831760044395,3.1812773627292903) (0.27367246913910304,3.0993097156366622) (0.2869741058408701,3.0177255303173585) (0.30018578312125094,2.9365219859654568) (0.3133100068148196,2.855696055705237) (0.32634923770193097,2.7752444732615555) (0.33930589577306575,2.695163703778193) (0.35218236424659677,2.6154499183985185) (0.36498099335399814,2.5360989722653953) (0.37770410390487813,2.4571063856361133) (0.39035399064273907,2.3784673278435555) (0.4029329254010089,2.3001766038669547) (0.41544316006766807,2.2222286433046374) (0.42788692936570133,2.1446174915673213) (0.44026645345562687,2.0673368031338706) (0.4525839403655061,1.9903798367322558) (0.46484158825309285,1.9137394523267897) (0.4770415875041622,1.837408109808764) (0.48918612267054007,1.76137786930152) (0.5012773742509489,1.68564039300286) (0.5133175203174813,1.6101869484976676) (0.5253087379903113,1.5350084134818676) (0.5372532047631412,1.460095281845401) (0.5491530996818753,1.3854376710670109) (0.5610106043790745,1.311025330877301) (0.5728279039669022,1.2368476531489878) (0.584607187791502,1.1628936829745582) (0.5963506500520378,1.0891521308918166) (0.608060490287986,1.0156113862172227) (0.6197389137386822,0.9422595314454951) (0.6313881315795677,0.869084357671953) (0.6430103610400872,0.7960733809914419) (0.6546078254086833,0.723213859824714) (0.6661827539308852,0.6504928131197822) (0.6777373816070168,0.5778970393722308) (0.689273948896586,0.5054131364048402) (0.7007947013369495,0.43302752184320575) (0.7123018890843412,0.3607264542204854) (0.7237977663858308,0.288496054640987) (0.7352845909912104,0.2163223289291741) (0.7467646235141959,0.14419119018781393) (0.7582401267526627,0.07208848168652934) (0.7697133649779095,0.0)};
\addlegendentry{Wrapped model}
\addplot[red!70!black,dashed,thick,domain=0:3.5,samples=100] {2*pi/1.2*exp(-x/sqrt(2))};
\addlegendentry{Gaussian parameter quotient}
\addplot[orange!85!black,dotted,thick,domain=.38:.769713364978] {2*pi*(.769713364978-x)};
\addlegendentry{Tangent at $\delta=D$, slope $-2\pi$}
\end{axis}
\end{tikzpicture}
\caption{Coordinate-circle lengths from $s_0=1.2$, with distance measured in each model. The wrapped circle closes at $D\approx0.7697$, with limiting slope $-2\pi$, while the quotient circle closes only at infinite distance.}\label{fig:cap}
\end{figure}

\FloatBarrier

\section{Discussion}
The two periodic constructions distinguish identifying parameters from changing observations. The Gaussian parameter quotient retains a complete hyperbolic metric. Wrapping strictly reduces information and moves the uniform limit to a finite intrinsic distance. Its completion is nevertheless regular: the uniform law is an analytic interior point of a two-dimensional statistical manifold.

The curvature gives a more precise distinction. Near the small-standard-deviation end the wrapped geometry approaches the Gaussian curvature $-1/2$, with differentiated exponentially small metric errors. Toward uniformity the curvature crosses zero once, reaches one positive maximum, and returns to zero at the smooth center. The last assertion concerns curvature at a point, not flatness on a neighborhood. The minimizing-diameter result also shows that adjoining this point restores minimizers for explicit pairs of nonuniform laws.

All statements here concern the full circular family, with mean and standard deviation varying independently. They use neither a higher-dimensional completion theorem nor a numerical certificate. Questions not settled here include explicit distances for arbitrary pairs and a global description of the cut locus. The analytic and geometric results above provide a starting point for those questions.

\appendix
\section{Proof of the circular curvature theorem}\label{app:curvature}
We prove Theorem~\ref{thm:unique}, using the Fisher coefficients in \eqref{eq:AB} and the curvature formula \eqref{eq:curvature}. The endpoint properties are given by Proposition~\ref{prop:curvature} and Corollary~\ref{cor:narrow}. A power-series estimate and a Poisson summation formula control overlapping intervals and together establish the global assertion.

\subsection{Fisher dissipation}
\begin{lemma}[Monotonicity of the angular Fisher coefficient]\label{lem:angular-monotonicity}
For $0<q<1$,
\begin{equation}\label{eq:Amonotone}
 qA'(q)=2\int_0^{2\pi}f\bigl((\log f)_{uu}\bigr)^2\dd u>0.
\end{equation}
\end{lemma}
\begin{proof}
Put $t=-\log q$, $v=(\log f)_u$, and $J(t)=\int fv^2\dd u=A(e^{-t})$.
The heat equation $f_t=f_{uu}$ gives $f_t=f(v_u+v^2)$ and
$v_t=v_{uu}+2vv_u$. Differentiation and periodic integration by parts yield
\[
 J'=\int f(2vv_{uu}+5v^2v_u+v^4)\dd u
    =-2\int fv_u^2\dd u,
\]
where $\int fvv_{uu}=-\int fv_u^2-\int fv^2v_u$ and
$\int(fv^3)_u=\int f(v^4+3v^2v_u)=0$.
On compact positive time intervals, normal convergence and a positive lower
bound for $f$ justify these operations. Since $J'=-qA'(q)$, the identity follows.
Equality would force $(\log f)_{uu}=0$; periodicity would then make $f$ constant,
contrary to its nonzero first Fourier moment $q$.
\end{proof}

\subsection{Lambert series for the Fisher coefficients}
\begin{lemma}[Lambert series for the Fisher coefficients]\label{lem:lambert}
Put $\xi=q^2$, $a(\xi)=A(q)/(2q^2)$, and $b(\xi)=B(q)/2$. Then
\begin{align}
 a(\xi)&=\sum_{n\ge1}(-1)^{n+1}
  \frac{n\xi^{n(n+1)/2-1}}{1-\xi^n},\label{eq:lambert-a}\\
 b(\xi)&=\sum_{n\ge1}(-1)^{n+1}
  \frac{n^2\xi^{n(n+1)/2-1}(1+\xi^n)}{(1-\xi^n)^2}.
 \label{eq:lambert-b}
\end{align}
These series are analytic for $|\xi|<1$. Their power-series coefficients are integers satisfying $|a_k|,|b_k|\le4(k+1)^2$. With primes now denoting $\xi$ derivatives,
\begin{equation}\label{eq:kappa}
 K(q)=\kappa(\xi)=
 \frac{a^2b'-3aa'b+\xi\{(a')^2b+aa'b'-2aa''b\}}{2a^2b^2}.
\end{equation}
In particular, this expression has no artificial singularity at $\xi=0$.
\end{lemma}
\begin{proof}
The Jacobi product formula \cite[Eq.~20.5.3]{DLMF} gives
\[
 \log(2\pi f)=c(q)+2\sum_{n\ge1}c_n(q)\cos(nu),\qquad
 c_n(q)=\frac{(-1)^{n+1}q^n}{n(1-q^{2n})}.
\]
All derivatives converge locally uniformly for $0<q<1$. Since $\E\cos(nu)=q^{n^2}$, the identity $A=-\E(\log f)_{uu}$ gives
$A=2\sum n^2c_nq^{n^2}$. Also $\E(\log f)_q=0$ implies
$c'+2\sum c_n'q^{n^2}=0$. Differentiating this equality and using
$B=-\E(\log f)_{qq}$ gives $B=2\sum n^2c_n'q^{n^2-1}$.
Substitution yields \eqref{eq:lambert-a}--\eqref{eq:lambert-b}.

Writing $k=n(n+1)/2-1+nm$, $m\ge0$, the contributions to $a_k,b_k$ are respectively
$(-1)^{n+1}n$ and $(-1)^{n+1}n^2(2m+1)$.
If $N(N+1)\le2(k+1)<(N+1)(N+2)$, then
$\sum_{n=1}^N n\le k+1$ and
$n^2(2m+1)=2(k+1)n-n^3\le2(k+1)n$.
Thus the stated common bound holds, and also proves normal convergence on every smaller complex disk. Finally substitute $A=2\xi a$, $B=2b$, and $\partial_q=2q\partial_\xi$ into \eqref{eq:curvature}.
\end{proof}

\subsection{Power-series estimates near the uniform law}
Put $N=a^2b'-3aa'b+\xi\{(a')^2b+aa'b'-2aa''b\}$ and
\begin{equation}\label{eq:Mnumerator}
 M=abN'-2(a'b+ab')N,\qquad \kappa'=\frac{M}{2a^3b^3}.
\end{equation}
The denominator is positive on $(0,1)$. Explicitly,
\begin{align}\label{eq:Mexpanded}
 M={}&a^3bb''+4a(a')^2b^2-5a^2a''b^2
 -2\xi a(a')^2bb'+3\xi a^2a''bb'+\xi a^2a'bb''\notag\\
 &-2\xi a^2a'''b^2+4a^2a'bb'-2\xi(a')^3b^2
 +4\xi aa'a''b^2-2a^3(b')^2-2\xi a^2a'(b')^2.
\end{align}

\begin{lemma}[A single critical radius for $q^2\le3/8$]\label{lem:local-analytic}
The function $M$ has exactly one zero $\xi_{\mathrm m}$ in $(0,3/8]$.
It lies in $(1/10,1/4)$, satisfies $M'(\xi_{\mathrm m})<0$, and
$M>0$ before this zero and $M<0$ after it throughout this interval.
\end{lemma}
\begin{proof}
The proof of Lemma~\ref{lem:lambert} gives the sharper coefficient bounds
$|a_k|\le k+1$ and $|b_k|\le2(k+1)^2$.
The first four coefficients are $(1,1,-1,1)$ and $(1,3,1,7)$, respectively.
Consequently the nonnegative-coefficient series
\[
 U=(1-\xi)^{-2}-\xi-2\xi^2-3\xi^3,\qquad
 V=2(1+\xi)(1-\xi)^{-3}-1-5\xi-17\xi^2-25\xi^3
\]
majorize their absolute coefficients. Define $R$ from the twelve terms in
\eqref{eq:Mexpanded} by replacing $a,b$ and their derivatives by $U,V$ and
by replacing every scalar coefficient by its absolute value.
If $M=\sum m_k\xi^k$ and $R=\sum R_k\xi^k$, then $|m_k|\le R_k$.
At $\xi=1/2$,
\[
 (U,U',U'',U''')=(21/8,43/4,83,750),\qquad
 (V,V',V'')=(105/8,477/4,1235),
\]
so
\begin{equation}\label{eq:Mnorm}
 \sum_{k\ge0}|m_k|2^{-k}\le R(1/2)
 =\frac{15064018893}{2048}<8\cdot10^6.
\end{equation}
The finite coefficients needed below are listed in Table~\ref{tab:Mcoeff}.
They follow from \eqref{eq:Mexpanded} and the coefficient formula in the
proof of Lemma~\ref{lem:lambert}: only its terms $n\le8$, expanded through
degree $42$, are needed. No sign assumption on subsequent coefficients is made.
\begin{table}[ht]
\centering\small
\caption{Integer coefficients $m_k$ of \eqref{eq:Mnumerator}.}\label{tab:Mcoeff}
\begin{tabular}{@{}rrrrrrrr@{}}\toprule
$k$&$m_k$&$k$&$m_k$&$k$&$m_k$&$k$&$m_k$\\\midrule
0&10&10&80892&20&3398640&30&167799360\\
1&20&11&-57912&21&-11595180&31&-93112088\\
2&-60&12&384444&22&10070386&32&-102790866\\
3&-640&13&-302712&23&-11635656&33&245596048\\
4&-3960&14&1223880&24&2677826&34&172611688\\
5&-9156&15&-1280496&25&5881652&35&-1138439880\\
6&-18276&16&2313534&26&-40112532&36&2054430634\\
7&-27648&17&-1629612&27&68094720&37&-1646381936\\
8&-16938&18&933840&28&-31812038&38&-210790014\\
9&-21952&19&-1818864&29&-61447560&39&635839280\\\bottomrule
\end{tabular}
\end{table}

For $0\le\xi\le1/10$, the table gives
$\sum_{k=1}^9|m_k|10^{-k}<4$, while \eqref{eq:Mnorm} gives
$\sum_{k\ge10}|m_k|\xi^k<8\cdot10^6/5^{10}<1$. Hence $M>5$ there.
For $1/10\le\xi\le1/4$, the first three terms of $M'/\xi^2$ are
$20/\xi^2-120/\xi-1920\le-1120$; this expression is decreasing on that interval.
Among $4\le k<24$ the positive terms have indices
$10,12,14,16,18,20,22$, and
\[
 \sum_{\substack{4\le k<24\\m_k>0}}km_k4^{-(k-3)}
 <50+18+5+1+1+1+1=77.
\]
Since $k2^{-k}$ decreases for $k\ge24$, the differentiated tail is bounded by
\[
 \sum_{k\ge24}k|m_k|\xi^{k-3}
 \le64\sum_{k\ge24}k2^{-k}|m_k|2^{-k}
 <\frac{64\cdot24\cdot8\cdot10^6}{2^{24}}<733.
\]
Thus $M'<-310\xi^2<0$. At $\xi=1/4$ the terms of degrees zero through four
sum to $-455/32<-14$. The positive terms of degrees five through $23$ sum
to less than $77/640<1$, and the remaining tail is less than
$8\cdot10^6/2^{24}<1/2$. Therefore $M(1/4)<0$.

It remains to control $[1/4,3/8]$. To bound its infinite tail more sharply, write
\begin{align*}
 U&=1+\xi+\xi^2+\xi^3+\frac{5\xi^4}{1-\xi}
                   +\frac{\xi^5}{(1-\xi)^2},\\
 V&=1+3\xi+\xi^2+7\xi^3+\frac{50\xi^4}{1-\xi}
          +\frac{22\xi^5}{(1-\xi)^2}+\frac{4\xi^6}{(1-\xi)^3}.
\end{align*}
Their derivatives and the definition of $R$ express $R$ as a nonnegative
linear combination of $\xi^j(1-\xi)^{-p}$ with $j\le28$ and $p\le15$.
For $k\ge40$ the polynomial terms contribute zero; for every other summand
the ratio of successive coefficients is
$1+(p-1)/(k+1-j)\le27/13$. With $r=3/8$ this proves
\[
 \sum_{k\ge40}|m_k|r^k\le\sum_{k\ge40}R_kr^k
 \le5R_{40}r^{40}<\frac14,
 \qquad R_{40}=4351246467529530.
\]
The last coefficient is obtained by expanding the displayed rational functions
in the twelve-term definition of $R$; the last inequality is an integer comparison.
The table also gives
\[
 \sum_{\substack{10\le k<24\\m_k>0}}m_kr^k<\frac{46}{5},\qquad
 |m_k|<2^k\ (24\le k<40),\qquad
 \sum_{k=24}^{39}|m_k|r^k<4(3/4)^{24}<\frac1{200}.
\]
The polynomial of degrees zero through four is decreasing for $\xi\ge1/4$,
and all coefficients of degrees five through nine are negative. Hence on
$[1/4,3/8]$,
$M<-14+46/5+1/200+1/4<0$.
The intermediate value theorem and the strict derivative bound on
$[1/10,1/4]$ now prove all the assertions.
\end{proof}

\subsection{An exact integral representation at the other end}
For this subsection set $u(\xi)=\xi a(\xi)$, $v(\xi)=\xi b(\xi)$,
$D=\xi\partial_\xi$, and
\[
 G(\xi)=1-24\sum_{n\ge1}\frac{n\xi^n}{1-\xi^n},\qquad
 P(y)=\prod_{n\ge1}(1-y^n)^3.
\]
The letter $u$ here denotes a function of $\xi$, not the angular integration variable.

\begin{lemma}[Differential identity for the Fisher coefficients]\label{lem:uv-identity}
For $0<\xi<1$,
\begin{equation}\label{eq:uv-identity}
 v=8Du-Gu+\tfrac14DG.
\end{equation}
\end{lemma}
\begin{proof}
Write $\ell=\log f$, $p=\ell_\theta$, and $w=\ell_{\theta\theta}$ at fixed $q$.
The theta product and half-period translation formulas
\cite[\S20.2(ii)--(iv), Eq.~20.2.13, and Eq.~20.5.1]{DLMF}
show that the meromorphic function $w$ is elliptic, with periods $2\pi$ and
$2\pi\tau$, where $q=e^{i\pi\tau}$. It has one double pole in a period cell.
At that pole, in a local coordinate $z$,
\[
 w=-z^{-2}-G/12+O(z^2).
\]
Indeed the product for $\theta_1(z/2,q)$ gives the constant
$-1/12+2\sum_{n\ge1}\xi^n/(1-\xi^n)^2=-G/12$.
Thus $w_{\theta\theta}+6w^2+Gw$ is an entire elliptic function and is constant.
The Fourier expansion
\[
 w=-2\sum_{n\ge1}\frac{(-1)^{n+1}nq^n}{1-q^{2n}}\cos(n\theta)
\]
shows that the unweighted period average of
$w_{\theta\theta}+6w^2+Gw$ is $-DG/2$. Consequently
\begin{equation}\label{eq:elliptic-w}
 w_{\theta\theta}+6w^2+Gw+\tfrac12DG=0.
\end{equation}
All series and derivatives used here converge normally on compact parameter
intervals and away from the complex poles.

Put $T=\E[w^2]$ and $Q=\E[p^4]$, with expectation under $f$.
Periodic integration by parts gives
\[
 \E w=-A,\qquad \E[p^2w]=-Q/3,\qquad
 \E[w_{\theta\theta}]=T-Q/3.
\]
The heat equation gives $q\partial_q f=-f_{\theta\theta}$, whence
$q^2B=T+Q/3$. The Fisher dissipation identity in
Lemma~\ref{lem:angular-monotonicity} gives $DA=T=2Du$.
Taking expectations in \eqref{eq:elliptic-w} yields
$7T-Q/3-GA+DG/2=0$. Substituting $A=2u$ and $q^2B=2v$
proves \eqref{eq:uv-identity}.
\end{proof}

\begin{lemma}[Dual integral formulas]\label{lem:dual-integral}
Let $h=s^2/\pi^2$, $\xi=e^{-\pi^2h}$, $y=e^{-4/h}$, and
\[
 I_j(h)=\int_{\R}z^je^{-hz^2/2}\operatorname{sech}z\dd z,\qquad
 \lambda(h)=\sqrt{\frac2{\pi h}}e^{-1/(2h)}.
\]
The normalized information deficits
$\eta=1-2\pi^2h\,u(\xi)$ and $\zeta=1-\pi^4h^2v(\xi)$ satisfy
\begin{align}
 \eta&=\lambda I_0 P(y)+\frac8h u(y),\label{eq:dual-eta}\\
 \zeta&=2h\lambda I_2 P(y)+\frac{16}h u(y)-\frac{16}{h^2}v(y).
 \label{eq:dual-zeta}
\end{align}
These identities hold for every $h>0$.
\end{lemma}
\begin{proof}
Set $T=\pi^2h$ and introduce the even Schwartz function
$g_T(z)=z e^{-Tz^2/2}/(2\sinh(Tz/2))$, with $g_T(0)=1/T$.
The Lambert series gives
$\eta=T\sum_{n\in\Z}(-1)^ng_T(n)$.
For the Fourier transform convention
$\widehat g(\omega)=\int e^{-i\omega z}g(z)\dd z$, the elementary transform
of $z/(2\sinh(Tz/2))$ is
$\pi^2T^{-2}\operatorname{sech}^2(\pi\omega/T)$.
Indeed, the substitution $t=e^{2z}$ and the beta integral give
\[
 \int_{\R}e^{-i\omega z}\operatorname{sech}^2z\dd z
 =2\Gamma(1-i\omega/2)\Gamma(1+i\omega/2)
 =\frac{\pi\omega}{\sinh(\pi\omega/2)},
\]
with the value at $\omega=0$ understood by continuity. Scaling and Fourier
inversion give the stated transform.
Convolving with the Gaussian transform and applying Poisson summation therefore gives
\begin{equation}\label{eq:poisson-deficit}
 \eta=\sqrt{\frac2{\pi h}}\sum_{k\ge0}e^{-(2k+1)^2/(2h)}
 \int_{\R}e^{-hz^2/2}\frac{\cosh((2k+1)z)}{\cosh^2z}\dd z.
\end{equation}
The elementary identity
\[
 \frac{\cosh((2k+1)z)}{\cosh^2z}
 =(-1)^k(2k+1)\operatorname{sech}z
 +4\sum_{j=1}^k(-1)^{k-j}(k-j+1)\cosh((2j-1)z)
\]
follows by induction from the addition formulas for $\cosh$.
The first term sums to $\lambda I_0P(y)$ by Jacobi's identity
$P(y)=\sum_{k\ge0}(-1)^k(2k+1)y^{k(k+1)/2}$, obtained by differentiating
the Fourier series and product for $\theta_1$ at zero
\cite[Eqs.~20.2.1 and 20.5.1]{DLMF}.
In the second term, perform the Gaussian integral and put $n=k-j+1$.
Its exponent becomes $y^{n(n+1)/2+n(k-n)}$, so this contribution is $8u(y)/h$.
The rearrangements are justified by
$\sum_{n\ge1}n y^{n(n+1)/2}/(1-y^n)<\infty$;
the corresponding $P$ series is absolutely convergent as well.
This proves \eqref{eq:dual-eta}.

The transformation of the Dedekind eta product \cite[Eq.~23.18.5]{DLMF},
or the theta transformation followed by differentiation at zero, gives
\[
 G(\xi)=\frac{12}{\pi^2h}-\frac4{\pi^2h^2}G(y),\qquad
 y\partial_y\log P(y)=\frac{G(y)-1}{8}.
\]
For clarity, the first formula also follows by logarithmically differentiating
\[
 \prod_{n\ge1}(1-e^{-Tn})
 =\sqrt{\frac{2\pi}{T}}\,
 e^{T/24-\pi^2/(6T)}\prod_{n\ge1}(1-e^{-4\pi^2n/T}).
\]
Using \eqref{eq:uv-identity}, these formulas imply
\[
 \zeta=(2G(y)/h-2)\eta-4h\partial_h\eta-4h^{-2}(D_yG)(y),
 \qquad D_y=y\partial_y.
\]
Insert \eqref{eq:dual-eta}, use $\partial_h I_0=-I_2/2$,
and apply \eqref{eq:uv-identity} at $y$. This gives \eqref{eq:dual-zeta}.
All differentiated series converge normally on compact subsets of $h>0$.
\end{proof}

\subsection{A uniform curvature inequality}
\begin{lemma}[Monotonicity up to $s=\pi/\sqrt{10}$]\label{lem:extended-endpoint}
For $0<s\le\pi/\sqrt{10}$ one has $s\,dK/ds>0$.
Equivalently, $\kappa'(\xi)<0$ for $e^{-\pi^2/10}\le\xi<1$.
\end{lemma}
\begin{proof}
Use $h,\eta,\zeta$ from Lemma~\ref{lem:dual-integral} and let
$\mathcal E=2h\partial_h=s\partial_s$. The metric is
$s^{-2}(1-\eta)\dd\alpha^2+2(1-\zeta)\dd(\log s)^2$.
With $U=1-\eta$, $V=1-\zeta$, and subscripts denoting $\mathcal E$ derivatives, set
\begin{align}\label{eq:curvature-loss-numerator}
 n&=(-2UU_2+4UU_1-4U^2+U_1^2)V+U(U_1-2U)V_1,\notag\\
 m&=(\mathcal En)UV-2n(U_1V+UV_1).
\end{align}
Direct substitution in the curvature formula gives
$\mathcal EK=m/(8U^3V^3)$, whose denominator is positive.

First replace the deficits by $\eta^{(0)}=\lambda I_0$,
$\zeta^{(0)}=2h\lambda I_2$, and denote the resulting numerator by $m^{(0)}$.
Put $\varepsilon=\lambda I_0$, $r=I_2/I_0$, $t=I_4/I_0$, $w=I_6/I_0$.
For $0<h\le1/10$,
\begin{equation}\label{eq:moment-box}
 0<\varepsilon<3/50,\qquad 0<r<5/2,\quad 0<t<31,\quad 0<w<920.
\end{equation}
Indeed each moment ratio decreases under the tilt $e^{-hz^2/2}$:
its derivative is minus one half of the covariance of the increasing functions
$|z|^{2j}$ and $|z|^2$, which is nonnegative by the two-independent-copies formula.
At $h=0$ the ratios are $\pi^2/4,5\pi^4/16,61\pi^6/64$, respectively.
These follow by differentiating
$\int e^{i\omega z}\operatorname{sech}z\dd z=\pi\operatorname{sech}(\pi\omega/2)$
at zero. The latter transform follows by integrating around a rectangle of height
$i\pi$. The bounds in \eqref{eq:moment-box} follow from $\pi<22/7$.
Also $\varepsilon\le\pi\lambda\le\sqrt{20\pi}e^{-5}<8/140<3/50$;
$\lambda$ increases on this interval and the exponential Taylor series gives $e^5>140$.

The seven derivatives divided by $\varepsilon$ are
\begin{align*}
 \eta^{(0)}/\varepsilon&=1,\\
 \eta^{(0)}_1/\varepsilon&=h^{-1}-1-hr,\\
 \eta^{(0)}_2/\varepsilon&=h^{-2}-4h^{-1}+1-2r+h^2t,\\
 \eta^{(0)}_3/\varepsilon&=h^{-3}-9h^{-2}+13h^{-1}-1+(-3h^{-1}+6-h)r\\
                  &\quad +(3h+3h^2)t-h^3w,\\
 \zeta^{(0)}/\varepsilon&=2hr,\\
 \zeta^{(0)}_1/\varepsilon&=2(1+h)r-2h^2t,\\
 \zeta^{(0)}_2/\varepsilon&=(2h^{-1}+2h)r-(4h+8h^2)t+2h^3w.
\end{align*}
Substitution in \eqref{eq:curvature-loss-numerator} defines polynomials $C_j$ by
$h^3m^{(0)}=\sum_{j=1}^5 \varepsilon^jC_j(h,r,t,w)$. The linear coefficient is
\[
 C_1=2\{1-11h+21h^2-3h^3+(-h^2+6h^3-3h^4)r
                   -(h^4+3h^5)t+h^6w\}>\frac9{50}.
\]
For this last inequality, $1-11h+21h^2-3h^3\ge107/1000$,
$h^2-6h^3+3h^4\le43/10000$, and
$h^4+3h^5\le13/100000$. These follow by differentiation on $[0,1/10]$;
thus the expression inside braces exceeds
$107/1000-(43/10000)(5/2)-(13/100000)31>9/100$.

The remaining finite polynomial estimates are recorded explicitly as follows.
If $Q_j$ is the sum of the absolute values of the negative-coefficient monomials
of $C_j/h$, substitution of $(h,r,t,w)=(1/10,5/2,31,920)$ gives
\begin{center}
\begin{tabular}{@{}ccc@{}}\toprule
$j$ & $Q_j(1/10,5/2,31,920)$ & Consequence\\\midrule
2 & $266563/20000<14$ & $C_2>6h$\\
3 & $85588293/2000000<43$ & $C_3>-43h$\\
4 & $19730389/2000000<10$ & $C_4>-10h$\\
5 & $5309537/8000000<1$ & $C_5>-h$\\\bottomrule
\end{tabular}
\end{center}
Here $C_2/h$ has the positive constant term $20$; all other positive monomials
were discarded. All powers in each $Q_j$ are nonnegative, so the displayed
single rational substitution bounds the whole box \eqref{eq:moment-box}.
The polynomials are specified by \eqref{eq:curvature-loss-numerator} and the
seven displayed derivatives, so these are finite algebraic inequalities,
with no parameter sampling. Since $6-43\varepsilon-10\varepsilon^2-\varepsilon^3>0$ for $\varepsilon<3/50$, it follows that
\begin{equation}\label{eq:leading-m-bound}
 h^3m^{(0)}>\frac9{50}\varepsilon.
\end{equation}

We next bound the exact corrections in \eqref{eq:dual-eta}--\eqref{eq:dual-zeta}.
For $0<y\le e^{-40}<1/100$ and $0\le j\le3$, the coefficient bounds already proved give
\[
 |D_y^ju(y)|\le2y,\qquad
 |D_y^j(P(y)-1)|\le6y,
 \qquad |D_y^jv(y)|\le4y\quad(0\le j\le2).
\]
For the product, use its Jacobi series and $2k+1\le3k(k+1)/2$ for $k\ge1$.
All these estimates reduce to
$\sum_{n\ge1}n^4y^n\le y/(1-16y)<2y$, since $n\le2^{n-1}$.
For any such function $F(y)$, with $|D_y^jF|\le By$,
\[
 \mathcal EF=\frac8hD_yF,\quad
 \mathcal E^2F=\frac{64}{h^2}D_y^2F-\frac{16}hD_yF,\quad
 \mathcal E^3F=\frac{512}{h^3}D_y^3F-\frac{384}{h^2}D_y^2F+\frac{32}hD_yF.
\]
Hence $h^j|\mathcal E^jF|\le600By$ for the applicable orders.
Leibniz's rule also gives
$h^j|\mathcal E^j(h^{-p}F)|\le1800B h^{-p}y$ for $p=1,2$ and $j\le3$.
The seven displayed derivatives and \eqref{eq:moment-box} give
$h^j|\eta_j^{(0)}|\le3\varepsilon$ for $j\le3$ and
$h^j|\zeta_j^{(0)}|\le \varepsilon$ for $j\le2$.
Applying the product rule to the exact formulas therefore yields
\begin{equation}\label{eq:scaled-dual-error}
 \max\left\{\max_{j\le3}h^j|\eta_j-\eta_j^{(0)}|,
             \max_{j\le2}h^j|\zeta_j-\zeta_j^{(0)}|\right\}
 \le2\cdot10^5h^{-2}e^{-4/h}=:\delta.
\end{equation}
For example, the product corrections are at most $86400\varepsilon y$ and $14400\varepsilon y$,
and the other corrections are at most $28800h^{-1}y$ and
$57600h^{-1}y+115200h^{-2}y$, respectively. These imply the stated bound.

As a polynomial in the seven loss variables, $m$ has degree at most five,
coefficient absolute sum $376$, and derivative weight at most three, assigning
weight $j$ to $\eta_j,\zeta_j$. These facts follow by expanding
\eqref{eq:curvature-loss-numerator}; its linear part is
$2\eta_3-4\eta_2+2\zeta_2-4\zeta_1$ and the remaining coefficient sum is $364$.
Thus, in the scaled variables $h^j\eta_j,h^j\zeta_j$, $h^3m$ has nonnegative
powers of $h$ and is $2000$-Lipschitz in the maximum norm on the unit box.
The segment between the approximate and exact scaled variables lies in that box:
its approximate entries have absolute value below $9/50$, and
$\delta\le2\cdot10^7e^{-40}<10^{-8}$.
Consequently
\[
 |h^3(m-m^{(0)})|\le2000\delta\le4\cdot10^8h^{-2}e^{-4/h}<\varepsilon/100.
\]
To verify the last inequality uniformly, $I_0\ge\pi-h\pi^3/8>2$, so $\varepsilon>2\lambda$.
The ratio of the penultimate expression to $\varepsilon$ is at most
$2\cdot10^8\sqrt{\pi/2}\,h^{-3/2}e^{-7/(2h)}$.
This increases on $(0,1/10]$ and is less than
$2\cdot10^8\cdot2\cdot40/10^{15}<1/100$ there,
using $e^{35}>(140)^7>10^{15}$.
Combining this with \eqref{eq:leading-m-bound} proves $m>0$ and hence
$\mathcal EK>0$. Finally $\mathcal E\xi=-2\pi^2h\xi<0$ gives the equivalent assertion.
\end{proof}

\begin{proof}[Proof of Theorem~\ref{thm:unique}]
Lemma~\ref{lem:local-analytic} proves that $\kappa'$ has precisely one zero on
$(0,3/8]$, is positive before it, and negative after it.
Lemma~\ref{lem:extended-endpoint} extends the negative sign to $1$.
The intervals overlap because $e^{-\pi^2/10}<3/8$:
$\pi>3.14$ and the first six terms of the exponential series give
$e^{\pi^2/10}>e^{197/200}>8/3$.
At the critical point, \eqref{eq:Mnumerator} gives
$\kappa''=M'/(2a^3b^3)<0$.
Since $\kappa(0)=0$, the maximum is positive; since
$\lim_{\xi\uparrow1}\kappa(\xi)=-1/2$, strict decrease after the maximum gives
exactly one zero, and its derivative is negative.
The change of variable $\xi=q^2$ preserves these signs and gives
$K''(q_{\mathrm m})=4q_{\mathrm m}^2\kappa''(q_{\mathrm m}^2)<0$.
\end{proof}

\section*{Use of generative AI}
During the preparation of this work, ChatGPT and Codex (OpenAI) were used to assist with proof exploration, drafting and revision, algebraic verification code, and \LaTeX{} preparation. The manuscript was subsequently reviewed and revised, and the proofs and calculations were checked. The authors take full responsibility for the content of this work.

\end{document}